\documentclass{ifacconf}

\usepackage{graphicx}      
\usepackage{natbib}        
\usepackage{amsmath,amssymb}
\usepackage{mathtools}

\newtheorem{remark}{Remark}

\begin{document}
\begin{frontmatter}

\title{The Innovation Null Space of the Kalman Predictor: A Stochastic Perspective for DeePC} 

\thanks[footnoteinfo]{This work was supported by VINNOVA Competence Center Dig-IT Lab Sustainable Industry, contract [2023-00556].}

\author[First]{Aihui Liu and Magnus Jansson} 

\address[First]{Department of Information Science and Engineering\\ 
KTH Royal Institute of Technology, 
   Stockholm, Sweden. \\
   (e-mail: aihui@kth.se; janssonm@kth.se).}

\begin{abstract}                
Willems' fundamental lemma uses a key decision variable $g$ to combine measured input-output data and describe trajectories of a linear time-invariant system. In this paper, we ask: what is a good choice for this vector $g$ when the system is affected by noise? For a linear system with Gaussian noise, we show that there exists an optimal subspace for this decision variable $g$, which is the null space of the innovation Hankel matrix. If the decision vector lies in this null space, the resulting predictor gets closer to the Kalman predictor. To show this, we use a result that we refer to as the Kalman Filter Fundamental Lemma (KFFL), which applies Willems’ lemma to the Kalman predictor. This viewpoint also explains several existing data-driven predictive control methods: regularized DeePC schemes act as soft versions of the innovation null-space constraint, instrumental-variable methods enforce it by construction, and ARX-based approaches explicitly estimate this innovation null space.
\end{abstract}

\begin{keyword}
Data-driven Control, DeePC, Willems' fundamental lemma
\end{keyword}

\end{frontmatter}


\section{Introduction}

Willems’ fundamental lemma is a deterministic statement about trajectories of a noise-free linear time-invariant (LTI) system: if the input is persistently exciting, then any finite-length trajectory of the system can be written as a linear combination of the columns of a Hankel matrix built from one sufficiently long, noise-free trajectory (\cite{Willems2005,Coulson2019}). In stochastic settings, however, the fundamental-lemma constraint becomes only approximate.

A common remedy for noisy data is to introduce regularization on the decision variables to avoid directions dominated by noise. Several works have studied the role of such \textbf{residual regularization} based on least squares residuals. These regularization terms can be formulated in various mathematically equivalent ways involving the variable $g$, see e.g. \cite{Breschi2023,bridging,Mattsson2024}. 

A second line of research draws on classic subspace identification, using \textbf{instrumental-variable (IV)} and oblique projections to cancel the effect of the innovations. Such ideas have been incorporated into data-driven predictive control to handle noise and closed-loop data, as in \cite{IV1,IV_Wang2023,IV3Dinkla2026}. 

A third line explicitly exploits the stochastic structure of the system via the innovation process of an optimal steady-state Kalman predictor: the \textbf{innovation-based} formulations in \cite{innoop} treat the innovation as an additional input and ensure that feasible $g$ generate zero future innovations. All these approaches implicitly restrict the subspace in which the decision variable $g$ is allowed to live, but a unifying view of how they choose this subspace is still missing.

This motivates the question: \textbf{what is an optimal subspace for the DeePC decision variable $g$ in the stochastic case?} We address this question by introducing the Kalman Filter Fundamental Lemma (KFFL). We apply Willems' fundamental lemma to the steady-state Kalman predictor, treating the measured outputs as additional inputs to a predictor system. Under Gaussian noise, we show that there exists a vector $g^*$ that satisfies the noisy fundamental lemma equations and reproduces the Kalman predictor. This $g^*$ lies in the null space of the innovation Hankel matrix $E_f$, which motivates interpreting $\text{null}(E_f)$ as a natural target subspace for $g$ in a stochastic DeePC setting. 
Informally, $\text{null}(E_f)$ collects those $g$ that are ``invisible'' to the innovations, and thus mimic what the Kalman predictor would do.


This perspective allows us to reinterpret several existing data-driven predictive control (DDPC) formulations as approximations or relaxations of the innovation null space condition. The main contributions of this paper are:

1. In Section \ref{sec:kffl}, we introduce the Kalman filter fundamental lemma (KFFL) and show that, under Gaussian noise, there exists a decision variable $g$ that reproduces the Kalman predictor and lies in $\text{null}(E_f)$, which we regard as an optimal subspace for $g$.

2. In Section \ref{sec:null}, we show that several DeePC and DeePC-like formulations can be understood as soft or hard relaxations of the innovation null space condition $E_fg=0$, thereby providing a unifying interpretation of regularization, instrumental-variable, and ARX-based approaches.

3. In Section \ref{sec:ARX}, we study ARX-based innovation estimation from this viewpoint, and discuss feasibility and some numerical issues.

We illustrate the idea of ARX-based nullspace estimation in a numerical example by analyzing the principal angles for the estimated subspace, and conclude with some final remarks.

\section{Problem setting}

\subsection{Notation}

For a signal sequence $\{ x_k \}_{k=0}^{N-1}$ of length $N$, we let $\mathbf{x}$ denote the stacked column vector $\mathbf{x}=[x_0^\top, \dots, x_{N-1}^\top]^\top$. For a sequence of matrices $X_1, \dots, X_n$, we denote $\text{col}(X_1,$ $\dots, X_n)$ $= [X_1^\top, \dots, X_N^\top]^\top$.

A signal trajectory $x_k \in \mathbb{R}^n$ is said to be persistently exciting of order $L$, if the block Hankel matrix 
\begin{equation*}
    X=\begin{bmatrix}
        x_0 & x_1 & \cdots & x_{N-L} \\
        x_1 & x_2 & \cdots & x_{N-L+1} \\
        \vdots & \vdots & \ddots & \vdots \\
        x_{L-1} & x_L & \cdots & x_{N-1}
    \end{bmatrix} \in \mathbb{R}^{nL \times (N-L+1)}
\end{equation*}
has full row rank. 

For a vector $x$, we define the $l_2$-norm as $||x||_2^2=x^\top x$, and the weighted $l_2$-norm as $||x||_W^2=x^\top W x$. For a matrix $M$, we let $\text{range}(M)$ denote its range space and $\text{null}(M)$ denote its null space. 

For a matrix $\Phi$ we denote the Moore-Penrose pseudoinverse as $\Phi^\dagger$, and define $\Pi = \Phi^\dagger\Phi$ and $\Pi_\perp = I-\Phi^\dagger \Phi$ as the orthogonal projectors onto the subspace spanned by rows of $\Phi$ and its orthogonal complement, respectively.

\subsection{Background on DeePC}

Consider a discrete-time LTI system subject to additive uncertainty:
\begin{equation} \label{eq:ss_form}
\begin{cases}
x_{t+1} = Ax_t + B u_t + w_t, \\
y_t = C x_t + Du_t+ v_t, 
\end{cases}
\end{equation}
where $x_t \in \mathbb{R}^{n},u_t \in \mathbb{R}^{n_u}$ and $y_t \in \mathbb{R}^{n_y}$ denote the state, input and output signals. $w_t \in \mathbb{R}^{n_x}$ and $v_t \in \mathbb{R}^{n_y}$ are the process and measurement noise with covariance matrices $\Sigma_w$ and $\Sigma_v$. 

Define Hankel matrices and the vectors
\begin{align} 
    U_p & = \begin{bmatrix}
        u_0 & u_1 & \cdots & u_{N-L_p} \\
        \vdots & \vdots & \ddots & \vdots  \\
        u_{L_p-1} & u_{L_p} & \cdots & u_{N-1}
    \end{bmatrix} \in \mathbb{R}^{n_uL_p \times (N-L_p+1)} \label{eq:hankel1} \\
    U_f &= \begin{bmatrix}
        u_{L_p}  & \cdots & u_{N} \\
        \vdots  & \ddots & \vdots \\
        u_{L_p+L_f-1}  & \cdots & u_{N+L_f-1}
    \end{bmatrix} \in \mathbb{R}^{n_uL_f \times (N-L_p+1)}  \label{eq:hankel2}
\end{align}
\begin{align}
    \mathbf{u}_{\text{ini}}(t) & = \text{col}(u_{t-L_p}, \dots, u_{t-1}) \in \mathbb{R}^{n_uL_p} \\
    \mathbf{u}(t) &= \text{col}(u_{t}, \dots, u_{L_f-1}) \in \mathbb{R}^{n_uL_f}
\end{align}
and similarly for $Y_p, Y_f, \mathbf{y}_{\text{ini}}$ and $\mathbf{y}$. 

For the predictive control problem, the aim is to minimize a finite-horizon control cost at each time step $t$, for example a quadratic cost of the form
\begin{equation} \label{eq:cost}
J(\mathbf{u},\hat{\mathbf{y}})=\Vert \hat{\mathbf{y}}-\mathbf{r}\Vert^2_Q + \Vert \mathbf{u} \Vert^2_R
\end{equation}
where $\mathbf{r}$ denotes a reference trajectory and $Q,R$ are user-defined positive definite weighting matrices. This optimization is performed under the input and output constraints $\mathbf{u} \in \mathcal{U}, \hat{\mathbf{y}} \in \mathcal{Y}$. 

Assume that the input sequence $u$ is persistently exciting so that the associated Hankel matrices have full row rank. For noise-free data, the result by Willems' fundamental lemma (\cite{Willems2005,Coulson2019}) implies that $(\mathbf{u}_{\text{ini}},\mathbf{u},\mathbf{y}_{\text{ini}},\hat{\mathbf{y}})$ is a trajectory of the noiseless system if and only if there exists a $g$ such that 
\begin{equation} \label{eq:deepc}
\begin{aligned}
    \text{col}(U_p, Y_p, U_f)g & =\text{col}(\mathbf{u}_{\text{ini}}, \mathbf{y}_{\text{ini}},\mathbf{u}) \\
    \hat{\mathbf{y}} &=Y_f g .
\end{aligned}
\end{equation}
Note that the size of $g$ grows with the number of training samples $N$. 

For noise-free trajectories, $g$ has no constraints other than (\ref{eq:deepc}). However, with noisy data, we need to further restrict $g$, and a regularization term is typically added to the control cost (\ref{eq:cost}). In Section \ref{sec:ls}, we will explore regularization methods such as those in \cite{bridging,Mattsson2023} and others.

\subsection{Steady-state Kalman filter system}

We can write the system (\ref{eq:ss_form}) in the innovation form:
\begin{equation} \label{eq:inno_form}
\begin{cases}
x_{t+1} = Ax_t + Bu_t + K e_t, \\
y_t = C x_t +D u_t+ e_t, 
\end{cases}
\end{equation}
with $e_t$ being the zero-mean Gaussian innovation.  

When the state measurements are not available, we can use a Kalman filter to estimate the state of the system. Under the assumption of Gaussian noise, i.e., $w(t)$ and $v(t)$ are normally distributed, the optimal state estimation in the mean-squared error sense is given by the steady-state Kalman filter (SSKF):
\begin{equation} \label{eq:sskf}
    \begin{aligned}
        \hat{x}_{t+1} &= A\hat{x}_t+Bu_t+K [y_t-\hat{y}_t] \\
        \hat{y}_t&= C\hat{x}_t+D u_t
    \end{aligned}
\end{equation}
where $K$ is the steady-state Kalman gain, which can be computed from the discrete-time algebraic Riccati equation(DARE), and all eigenvalues of $A-KC$ lie strictly inside the unit circle under weak conditions.



\section{Kalman filter fundamental lemma} \label{sec:kffl}

We now apply Willems’ fundamental lemma to the steady-state Kalman filter viewed as a predictor system. We rewrite the SSKF system (\ref{eq:sskf}) by treating the measured output as an additional input and the predicted output as the output:
\begin{equation} \label{eq:kfflsys}
    \begin{aligned}
        \hat{x}_{t+1}&= (A-KC)\hat{x}_t + \begin{bmatrix}
            B-KD & K
        \end{bmatrix} \begin{bmatrix}
            u_t \\ y_t
        \end{bmatrix} \\
        \hat{y}_t & = C\hat{x}_t + \begin{bmatrix}
            D & 0
        \end{bmatrix} \begin{bmatrix}
            u_t \\ y_t
        \end{bmatrix}
    \end{aligned}
\end{equation}
We regard $(u,y)$ as the inputs and $\hat{y}$ as the output of this predictor system, and build the Hankel matrices $\hat{Y}_p,\hat{Y}_f, E_p,E_f$ accordingly. 

\begin{lem} \label{lemma} (Kalman Filter Fundamental lemma)
Under the rank and PE assumptions, we apply Willems’ fundamental lemma to the SSKF system (\ref{eq:kfflsys}). For any trajectory $(\mathbf{u}_{\text{ini}},\mathbf{y}_{\text{ini}},\hat{\mathbf{y}}_{\text{ini}}, \mathbf{u}, \mathbf{y}, \hat{\mathbf{y}})$ of (\ref{eq:kfflsys}), there exists $g$ such that
    \begin{equation} 
    \text{col}(U_p, Y_p,\hat{Y}_p,
        U_f,Y_f,\hat{Y}_f)g = \text{col}(\mathbf{u}_{\text{ini}},\mathbf{y}_{\text{ini}},\hat{\mathbf{y}}_{\text{ini}}, \mathbf{u}, \mathbf{y}, \hat{\mathbf{y}} ) .
\end{equation}
\end{lem}
Since $y$ is noisy, $Y_p$ and $Y_f$ have full row rank almost surely, so $\text{col}(U,Y)$ is persistently exciting as soon as $u$ is. 

Next, we view this trajectory constraint as a predictor. For prediction, we consider hypothetical future trajectories with zero measurement noise, so that the predicted output $\hat{\mathbf{y}}$ coincides with the true future output, i.e., $\mathbf{y}=\hat{\mathbf{y}}$ for the multi-step ahead predictor. In that case, Lemma \ref{lemma} specializes to 
\begin{equation} \label{eq:kffl}
    \text{col}(U_p, Y_p,\hat{Y}_p,
        U_f,Y_f,\hat{Y}_f)g = \text{col}(\mathbf{u}_{\text{ini}},\mathbf{y}_{\text{ini}},\hat{\mathbf{y}}_{\text{ini}}, \mathbf{u}, \hat{\mathbf{y}}, \hat{\mathbf{y}} ) .
\end{equation}
Subtracting the last two constraints, we get $(Y_f-\hat{Y}_f)g=E_fg=0$, where $E_f$ is the future innovation Hankel matrix. Thus, any $g$ that generates a valid trajectory for the predictor with zero future innovation must satisfy $E_f g = 0$. This is the innovation null-space condition.

\begin{remark} (Connection to the innovation-based formulation).
    The idea in Lemma \ref{lemma} was originally given in \cite{innoop}. They apply Willems' fundamental lemma to the innovation system (\ref{eq:inno_form}) by viewing the innovation $e_t$ as additional input, and by setting the future innovations $\mathbf{e}=0$:
    \begin{equation} \label{eq:innoKF}
        \text{col}(U_p,Y_p,E_p,U_f,Y_f,E_f)g=\text{col}(\mathbf{u}_{\text{ini}}, \mathbf{y}_{\text{ini}}, \mathbf{e}_{\text{ini}}, \mathbf{u}, \hat{\mathbf{y}}, 0) ,
    \end{equation}
    which is equivalent to (\ref{eq:kffl}). Their last constraint also leads to $E_fg=0$. 
\end{remark}

Our goal is to identify the subspace that captures the ``good'' choices of the decision variable $g$ in the stochastic setting. The following is a Kalman-predictor specialization of the noisy Willems’ fundamental lemma. We refer to this as the Kalman Filter Fundamental Lemma (KFFL).

\begin{thm} \label{thm1}
(Innovation null space constraint in noisy DeePC) 
Consider the LTI system (\ref{eq:ss_form}) with independent, zero-mean Gaussian process and measurement noise, and let $E_f$ denote the future innovation Hankel matrix. For any given $(\mathbf{u}_{\text{ini}},\mathbf{u},\mathbf{y}_{\text{ini}})$, let $\hat{\mathbf{y}}_{\text{KF}}$ be the future output sequence predicted by the steady-state Kalman filter. Then there exists a vector $g^*$ such that
\begin{equation} 
\begin{aligned}
    \text{col}(U_p, Y_p, U_f)g^* & =\text{col}(\mathbf{u}_{\text{ini}}, \mathbf{y}_{\text{ini}},\mathbf{u}) \\
    \hat{\mathbf{y}}_{\text{KF}} &=Y_f g^* ,
\end{aligned}
\end{equation}
and this $g^*$ lies in the null space of the innovation Hankel matrix:
    \begin{equation} \label{Ef}
        g^* \in \text{null}(E_f) .
    \end{equation}
\end{thm}

\begin{pf}
The SSKF system (\ref{eq:kfflsys}) is an LTI system with inputs $(u_t, y_t)$ and output $\hat{y}_t$. Under the persistency-of-excitation conditions of Willems' fundamental lemma, Lemma \ref{lemma} guarantees that for any trajectory $(\mathbf{u}_{\text{ini}},\mathbf{y}_{\text{ini}},\hat{\mathbf{y}}_{\text{ini}}, $ $ \mathbf{u}, \hat{\mathbf{y}}, \hat{\mathbf{y}})$ of (\ref{eq:kfflsys}), there exists a $g$ such that (\ref{eq:kffl}) holds. By partitioning (\ref{eq:kffl}) into past-future blocks and separating measured and predicted outputs, we get
\begin{subequations} \label{eq:kffl_sep}
    \begin{align}
        \text{col}(U_p, Y_p, U_f)g &= \text{col}(\mathbf{u}_{\text{ini}},\mathbf{y}_{\text{ini}}, \mathbf{u} )  \label{eq:kffl_1} \\
    \hat{Y}_pg&= \hat{\mathbf{y}}_{\text{ini}} \label{eq:kffl_2} \\
    (Y_f-\hat{Y}_f)g &=0 \label{eq:kffl_3} \\
    \hat{\mathbf{y}}_{\text{KF}} &= Y_f g . \label{eq:kffl_4}
    \end{align}
\end{subequations} 
Observe that the noisy DeePC constraints (\ref{eq:deepc}) are contained in (\ref{eq:kffl_1}) and (\ref{eq:kffl_4}). Hence, Lemma \ref{lemma} guarantees the existence of a Kalman filter optimal $g^*$ such that (\ref{eq:deepc}) holds with $\hat{\mathbf{y}} = \hat{\mathbf{y}}_\text{KF}$. For this vector $g^*$, (\ref{eq:kffl_3}) implies $(Y_f-\hat{Y}_f)g^* =E_fg^*=0$. This finishes the proof of $g^* \in \text{null}(E_f)$. 
\end{pf}

\begin{remark}
    For a given predictor system, we can view $u_t$ and $y_t$ as inputs, and the estimated $\hat{y}_t$ as the output of this predictor system. The same construction applies to any predictor (such as ARX), which we exploit in Section \ref{sec:null}.
\end{remark}

In summary, within the set of $g$ that satisfy the noisy fundamental lemma equations (\ref{eq:deepc}), one can always find a $g$ that reproduces the Kalman predictor’s output and is ``invisible'' to the innovations. Since the Kalman predictor is the mean-square optimal predictor under Gaussian noise, we interpret $\text{null}(E_f)$ as an optimal subspace for $g$.


\section{the innovation null space $E_f$} \label{sec:null}

In this section, we survey three ways in which existing DDPC schemes relate to the innovation null space $\text{null}(E_f)$: least-squares residual soft regularization, instrumental-variable formulations, and ARX-based innovation hard constraint. 


\subsection{Estimate $\text{null}(E_f)$ via least-squares residuals} \label{sec:ls}

In this subsection, we approximate the innovation Hankel matrix $E_f$ by an ordinary least-squares (OLS) residual, and then enforce the innovation null space condition $E_{f} g = 0$ from Theorem \ref{thm1} in a relaxed way via regularization.

Guided by Theorem \ref{thm1}, our task in the stochastic setting is to estimate the innovation null space $\text{null}(E_f)$ from data, and restrict $g$ to be close to this subspace. In this subsection, we focus on open-loop data, where the innovation $E_f$ is uncorrelated with the inputs. We show that, in this case, the directions in $Y_f$ associated with the innovations can be approximated by the directions that are orthogonal to the data subspace spanned by
\begin{equation} \label{eq:Phi}
    \Phi = \text{col}(U_p , Y_p , U_f) .
\end{equation}
Intuitively, the innovation lives in the part of $Y_f$ that the regressor $\Phi$ cannot explain. Several existing DeePC-type formulations can be interpreted as constraining $g$ so that it has only a small component in this orthogonal complement. 

We will first construct a residual-based approximation of $E_f$, then discuss how to enforce the innovation null space condition softly via regularization.

\subsubsection*{4.1.1 Residual-based approximation of  $E_f$} 

We start from the subspace predictive control (SPC) formulation in \cite{SPC}. The SPC predictor identifies a linear mapping $\Theta$ from the regressor $\Phi$ to the future output Hankel matrix $Y_f$ by solving the least squares problem
\begin{equation}
    \min_{\Theta} \Vert  (Y_f - \Theta \Phi ) \Vert^2 .
\end{equation}
The ordinary least squares solution is
\begin{equation}
    \hat{\Theta}_{\text{LS}} = Y_f \Phi^\dagger
\end{equation}
and the corresponding estimated future-output Hankel matrix is
\begin{equation}
\hat{Y}_f=\hat{\Theta}_{\text{LS}}\Phi=Y_f\Phi^\dagger \Phi=Y_f \Pi ,
\end{equation}
where $\Pi=\Phi^\dagger \Phi$ is the orthogonal projector onto $\text{range}(\Phi)$. 

Given $(\mathbf{u}_{\text{ini}},\mathbf{y}_{\text{ini}}, \mathbf{u})$, the SPC predictor is
\begin{equation}
\hat{\mathbf{y}}_{\text{LS}}=\hat{\Theta}_{\text{LS}} \text{col}(\mathbf{u}_{\text{ini}},\mathbf{y}_{\text{ini}}, \mathbf{u} ) .
\end{equation}
The associated residual Hankel matrix is
\begin{equation}
    \hat{E}_{f,\text{LS}} = Y_f(I-\Phi^\dagger \Phi)=Y_f(I-\Pi) .
\end{equation}
We interpret $\hat{E}_{f,\text{LS}}$ as a residual-based approximation of the innovation Hankel matrix $E_f$, obtained as the part of $Y_f$ that is orthogonal to the regressor space spanned by$\Phi$. 

Using $\hat{E}_{f,\text{LS}}$ as a proxy for $E_f$ in the innovation null space condition of Theorem \ref{thm1} suggests the approximate constraint 
\begin{equation}
    \hat{E}_{f,\text{LS}}g=Y_f(I-\Pi)g \approx 0 ,
\end{equation}
which encourages $g$ to lie in $\text{null}(\hat{E}_{f,\text{LS}})$, i.e., in an estimated version of $\text{null}(E_f)$. Under standard OLS assumptions for open-loop data, the LS residuals are consistent estimators of the innovation sequence. It is therefore natural to view $\text{null}(\hat{E}_{f,\text{LS}})$ as a proxy for the true innovation null space $\text{null}(E_f)$.





In summary, we know that the LS residual-based approximation suggests the component $Y_f(I-\Pi)g$ should be small. This can lead to a category of regularization choices.

\subsubsection*{4.1.2 Orthogonal-projection regularization methods} 

Next, we revisit several DDPC schemes and interpret them through this LS nullspace lens.

The first example is the orthogonal projection regularizer used in \cite{bridging}, which adds a term of the form $\lambda \Vert (I-\Pi)g \Vert_p$. Regardless of the choice of $p$-norm, this regularizer constrains $(I-\Pi)g$ to be small, and suppresses the ``innovation-like'' directions of $g$. It hence encourages $g$ to stay close to the approximated innovation null space $\hat{E}_{f,\text{LS}}g \approx 0$. 

\paragraph*{Splitting $g$.}

\cite{Mattsson2023,Mattsson2024} decompose $g$ into a least-squares prediction part and an innovation part by writing 
\begin{equation}
    g=\hat{g}_{\text{LS}}+g_\perp=\Phi^\dagger \text{col}(\mathbf{u}_{\text{ini}},\mathbf{y}_{\text{ini}}, \mathbf{u} )+ (I-\Pi)g'
\end{equation}
where $\hat{g}_{\text{LS}}(\mathbf{u}) \in \text{range}(\Phi^\top)$ is the SPC/OLS predictor, and $g_\perp=(I-\Pi)g'$ lies in the orthogonal complement $\text{null}(\Phi)$, with $g'$ a free coordinate vector parameterizing this innovation component.

Using this decomposition, the prediction becomes 
\begin{equation}
\begin{aligned}
    \hat{\mathbf{y}}=Y_fg & = \hat{\Theta}_{\text{LS}}\text{col}(\mathbf{u}_{\text{ini}},\mathbf{y}_{\text{ini}}, \mathbf{u} ) + \hat{E}_{f,\text{LS}}g' \\
    & = \hat{\mathbf{y}}_{\text{LS}} + \Delta \hat{\mathbf{y}}
\end{aligned}
\end{equation}
with $\Delta \hat{\mathbf{y}} = \hat{E}_{f,\text{LS}} g' \in \text{range}(\hat{E}_{f,\text{LS}})$. Thus, all deviations from the SPC/OLS predictor $\hat{\mathbf{y}}_{\text{LS}}$ that are subject to regularization are generated by $\hat{E}_{f,\text{LS}}$ and live in the residual-based innovation space. 

The DDPC formulation in \cite{Mattsson2024} optimizes over $(\mathbf{u},\Delta \hat{\mathbf{y}})$ by the cost function
\begin{equation} \label{eq:split}
J(\mathbf{u},\hat{\mathbf{y}})+\lambda_1 \Vert \hat{g}_{\text{LS}}(\mathbf{u}) \Vert^2 + \lambda_2 \Vert \hat{E}_{f,\text{LS}} g' \Vert^2_{ (\hat{E}_{f,\text{LS}}\hat{E}_{f,\text{LS}}^\top)^\dagger }
\end{equation}
subject to $\hat{\mathbf{y}}= \hat{\mathbf{y}}_{\text{LS}} +\hat{E}_{f,\text{LS}} g'$. This can be interpreted as weighted regularization, acting only along directions in $\text{range}(\hat{E}_{f,\text{LS}})$.


\paragraph*{$\gamma$-DDPC.}

The $\gamma$-DDPC formulation of \cite{Breschi2023cdc,Breschi2023} uses an LQ-decomposition of the stacked Hankel matrix
\begin{equation}
    \begin{bmatrix} \begin{bmatrix}
        U_p\\ Y_p
    \end{bmatrix} \\ U_f \\ Y_f \end{bmatrix}= \underbrace{ \begin{bmatrix}
L_{11} & 0      & 0      \\
L_{21} & L_{22} & 0      \\
L_{31} & L_{32} & L_{33} \end{bmatrix} }_{L}
    \underbrace{ \begin{bmatrix}
Q_{1} \\
Q_{2} \\
Q_{3} \end{bmatrix} }_{Q} 
\end{equation}
where $L$ is lower lock-triangular and $Q$ is orthogonal. Defining $\gamma = \text{col}(\gamma_1,\gamma_2,\gamma_3)  \coloneq  Qg$ leads to an optimization problem directly in terms of $\gamma$ with cost 
\begin{equation}
    J(\mathbf{u},\hat{\mathbf{y}})+ \beta_2 \Vert \gamma_2 \Vert^2+ \beta_3 \Vert \gamma_3 \Vert^2
\end{equation}
and linear constraints linking $(\mathbf{u},\hat{\mathbf{y}})$ to $\gamma$.
 
It can be shown (see \cite{Breschi2023cdc} and \cite{Mattsson2024}, or Appendix \ref{sec:appA}) that, for the optimal input $\mathbf{u}$, the prediction can be written as 
\begin{equation} \label{eq:gamma_de}
    \hat{\mathbf{y}} =\hat{\mathbf{y}}_{\text{LS}}+L_{33} \gamma_3 ,
\end{equation}
and that $L_{33}L_{33}^\top = Y_f(I-\Pi)Y_f^\top = \hat{E}_{f,\text{LS}}\hat{E}_{f,\text{LS}}^\top$. This implies that $\gamma_3$ parametrizes exactly the component of $\hat{\mathbf{y}}$ that lies in the residual-based innovation space $\text{range}(\hat{E}_{f,\text{LS}})$. Consequently, the regularization term $\beta_3 \Vert \gamma_3 \Vert^2$ enforces the solution to stay close to $\text{null}(\hat{E}_{f,\text{LS}})$ in a soft sense.

In summary, in this Subsection \ref{sec:ls}, the regularization terms in several different formulations tend to drive $\hat{E}_{f,\text{LS}}g$ towards zero. This yields a soft enforcement of $E_fg \approx 0$, and therefore a relaxation of the innovation null space condition in Theorem \ref{thm1}.

\subsection{Estimate $\text{null}(E_f)$ via instrumental variables} \label{sec:iv}

An alternative, rooted in classical subspace identification, is to use oblique projections and instrumental variables (IVs) to construct a subspace in which the innovation does not appear. \cite{IV1} introduce the idea of DeePC with an instrumental matrix $Z$. The IV condition is summarized as 
\begin{equation} \label{eq:instru}
    \lim_{N\to \infty} \frac{1}{N} E_f Z^\top \to 0, \quad \text{rank}(\Phi Z^\top)=\text{rank}(\Phi)
\end{equation}
so that $Z$ is asymptotically uncorrelated with the future innovations $E_f$. 

Geometrically, the first condition in (\ref{eq:instru}) states that the empirical cross-covariance between $E_f$ and $Z^\top$ tends to zero as $N\to \infty$. This implies 
\begin{equation}
    \text{range}(Z^\top) \subseteq \text{null}(E_f) 
\end{equation}
and hence any vector of the form $g=Z^\top h$ automatically satisfies the nullspace constraint $E_fg=E_f Z^\top h=0$.

The DeePC constraint (\ref{eq:deepc}) with IV becomes:
\begin{equation}
\begin{aligned}
    \text{col}(U_p, Y_p, U_f)Z^\top h & =\text{col}(\mathbf{u}_{\text{ini}}, \mathbf{y}_{\text{ini}},\mathbf{u}) \\
    \hat{\mathbf{y}} &=Y_f Z^\top h .
\end{aligned}
\end{equation}
Such instrumental variables realize the nullspace constraint $E_fg=0$ by construction, instead of by adding regularizations and penalties.

Specifically for open-loop data, \cite{IV1} uses the instrumental matrix $Z$ to be the SPC regressor $Z=\Phi$ in (\ref{eq:Phi}). However, this is not feasible in a closed-loop setting because of the correlation between the input and the innovation (\cite{IV2}). To handle closed-loop data, \cite{IV_Wang2023} propose two IV choices: the future reference signal and the controller's left coprime factorization (LCF). Similarly, \cite{IV3Dinkla2026} address the closed-loop setting by constructing $Z$ from suitably time-shifted past data.

\subsection{Estimate $\text{null}(E_f)$ via ARX residual} \label{sec:ARX}

The third class of methods estimates the innovation Hankel $E_f$ explicitly by fitting a high-order one-step-ahead ARX model. This follows the innovation-based approach in \cite{innoop}, where the innovation is treated as an additional input and Willems’ fundamental lemma is applied to the augmented system. The corresponding predictor is equivalent to applying our fundamental-lemma viewpoint to the innovation system and leads to the same innovation null space constraint ${E}_{f} g=0$. For simplicity, in this section, we write $\hat{E}_f$ and $\hat{Y}_f$ for the corresponding Hankel matrices from the one-step-ahead prediction, and $\rho$ for the ARX model order.


To exploit the constraint $\hat{E}_{f} g=0$, \cite{innoop} parametrizes $g$ via an order reduction. Let $\hat{E}_f^\perp$ be a basis matrix for the null space of $\hat{E}_f$, so that every feasible $g$ satisfying $\hat{E}_fg=0$ can be written as $g=\hat{E}_f^\perp h$ for some reduced-dimension coordinate $h$. Substituting this parametrization into the augmented fundamental lemma system yields
\begin{equation} 
    \text{col}(
        U_p \hat{E}_f^\perp , Y_p \hat{E}_f^\perp , \hat{E}_p \hat{E}_f^\perp ,
        U_f \hat{E}_f^\perp ) h = \text{col}(
        \mathbf{u}_{\text{ini}} , \mathbf{y}_{\text{ini}} ,\hat{\mathbf{e}}_{\text{ini}} ,
        \mathbf{u}  )
\end{equation}
which can be solved for $h_{\text{pinv}}$ by least squares or pseudo-inverse. The corresponding output prediction is then $\hat{\mathbf{y}}= Y_f \hat{E}_f^\perp h_{\text{pinv}}$.



From the KFFL viewpoint, this is a null space method: if $\text{null}(\hat{E}_f)$ matches the true innovation null space $\text{null}(E_f)$, the resulting predictor reproduces the Kalman filter optimal prediction as in Section \ref{sec:kffl}. In that idealized case, the ARX-based approach can be interpreted as an implementation of the KFFL constraint using an explicit parametrization of $\text{null}(E_f)$.

\begin{remark}
    The ARX-based predictor itself satisfies a fundamental-lemma–type data equation, so (\ref{eq:kffl}) is automatically feasible if both $\hat{\mathbf{e}}_{\text{ini}}(0)$ and $\hat{E}_p$ are obtained from the same ARX model. The Hankel matrices built from the corresponding ARX model have a row dimension that scales with the ARX order $\rho$. To satisfy the persistency-of-excitation conditions in \cite{innoop}, it is natural to choose $\rho \geq L_p$.
\end{remark}

In practice, however, $\text{null}(\hat{E}_f)$ is obtained from an estimated ARX model, and its null space depends on the model order $\rho$ and on finite-data effects. This highlights an important subtlety: although high-order ARX models yield consistent one-step-ahead predictions with enough data, the corresponding null space estimate $\text{null}(\hat{E}_f)$ is not guaranteed to be consistent with the true innovation null space. In the experiments, we will illustrate this effect by analyzing the principal angles between the estimated and true null spaces.

\section{Experiments} \label{sec:ex}

In this section, numerical simulations are carried out to show the estimation qualities of $\text{null}(\hat{E}_f)$ using various methods. We consider the LTI system presented in \cite{innoop} with the system matrices
\begin{equation}
    \begin{aligned}
        A& = \begin{pmatrix}
            0.7326 & -0.0861     \\
            0.1722 & 0.9909
        \end{pmatrix}, \quad B  = \begin{pmatrix}
        0.0609 \\ 0.0064 \end{pmatrix} , \\
        C & = \begin{pmatrix}
            0 & 1.4142 \end{pmatrix}, \quad D=0 ,
    \end{aligned}
\end{equation}
the process noise variance $\Sigma_w = \sigma_w^2 \mathbb{I}_2$, and the measurement noise variance $\Sigma_v = \sigma_v^2$. We set $\sigma_w=5 \times 10 ^{-3}$ and $\sigma_v=2\times 10 ^{-3}$. 

Firstly, we show that when using ARX to estimate the innovation null space $\text{null}(\hat{E}_{f})$, the estimation does not consistently improve as the ARX order $\rho$ increases. 

We collect closed-loop training trajectories of length $N_{\text{train}} = 200$ using the simple feedback $u(t) = 5(r_{\text{train}}(t)-y(t))$, where the external input $r_{\text{train}}(t)$ is a square wave with a period of $50$ time-steps and amplitude of $2$, contaminated by a zero-mean Gaussian distributed sequence with variance $0.01$. The past and future horizons are $L_p=10$ and $L_f=15$, and the control cost metrics $Q = 1$ and $R = 0.01$. We want the controlled output to follow a sinusoid reference signal $r(t)=\sin (2\pi t/N_{\text{test}})$ where $N_{\text{test}} = 100$.

\begin{figure}[tb]
    \centering
    \includegraphics[width= \linewidth]{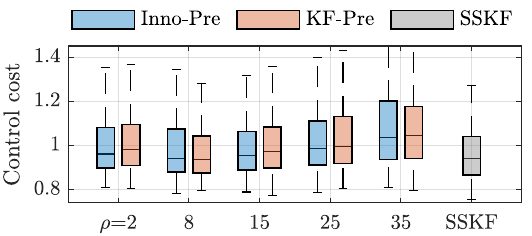}\\
    \includegraphics[width= \linewidth]{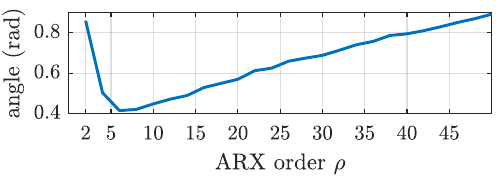}
    \caption{(Top) Total control cost as a function of increasing ARX order $\rho$. (Bottom) Largest principal angle between real $\text{null}(E_f)$ and ARX estimation $\text{null}(\hat{E}_{f})$.  Training data length $N_{\text{train}} = 200$ and Monte Carlo simulations $N_{MC}=300$.}
    \label{fig:fig1}
\end{figure}


In the top plot of Fig.~\ref{fig:fig1}, we compare the closed-loop control cost of three controllers: the controller \textbf{Inno-Pre} with an innovation predictor built from (\ref{eq:innoKF}), its equivalent variant \textbf{KF-Pre} using the predictor built from (\ref{eq:kffl}), and the Kalman filter baseline \textbf{SSKF}. The top plot shows that there exists an ARX order $\rho$ that yields the minimum control cost, and when we further increase $\rho$, the cost also increases. In the bottom plot of Fig.~\ref{fig:fig1}, we show the subspace estimation quality using the largest principal angle between the true innovation null space $\text{null}(E_f)$ and the estimated null space $\text{null}(\hat{E}_{f})$. We observe a similar trend: there exists a good choice of $\rho$ that gives the best $\text{null}(\hat{E}_{f})$ estimation, and further increasing the ARX order will not help.

\begin{figure}[tb]
\begin{center}
\includegraphics[width=8.4cm]{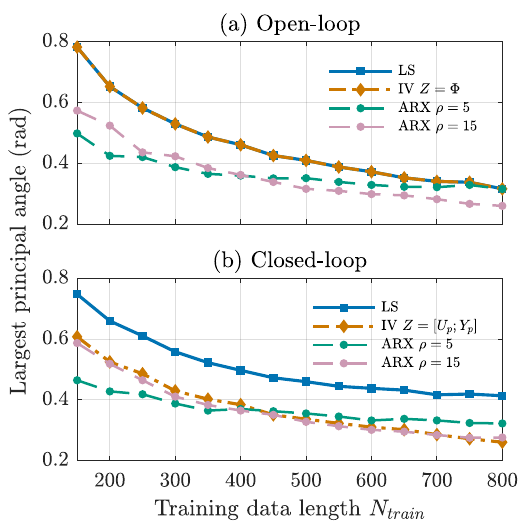}    
\caption{Largest principal angle between real $\text{null}(E_f)$ and estimated $\text{null}(\hat{E}_f)$ for different training data length $N_{\text{train}}$, for (a) open-loop and (b) closed-loop training data. The LS estimation $\text{null}(\hat{E}_{f,\text{LS}})$ is biased for closed-loop data.} 
\label{fig:fig2}
\end{center}
\end{figure}

Then in Fig.~\ref{fig:fig2}, we compare the largest principal angle between the true null space $\text{null}(E_f)$ and the estimated $\text{null}(\hat{E}_{f})$ obtained by different methods: least squares residual-based \textbf{LS} from Section \ref{sec:ls}, instrumental variables \textbf{IV} from Section \ref{sec:iv}, and ARX-based estimation \textbf{ARX} from Section \ref{sec:ARX}. The open-loop training data in Fig.~\ref{fig:fig2}(a) is generated by the input $u_{OL}(t) = r_{\text{train}}(t)$, and the closed-loop training data in Fig.~\ref{fig:fig2}(b) by the input $u_{CL}(t) = 5(r_{\text{train}}(t)-y(t))$. For IV options, we use $Z=\Phi$ for open-loop data, and $Z=\text{col}(U_p,Y_p)$ for closed-loop data. 

From Fig.~\ref{fig:fig2}, we observe that as we have more training data, the largest principal angle between $\text{null}(E_f)$ and the estimated $\text{null}(\hat{E}_f)$ decreases, i.e., the subspace estimation gets closer to the real $\text{null}(E_f)$. For open-loop data, using the least-squares residual is equivalent to the IV methods when the IV $Z=\Phi$ is the LS regressor. However, for closed-loop data, the LS residual is biased away from the innovation, and the estimation $\text{null}(\hat{E}_{f,\text{LS}})$ is also biased.

\section{Conclusion} \label{sec:conclu}

In this paper we introduced the Kalman Filter Fundamental Lemma (KFFL), showing that the innovation null space $\text{null}(E_f)$ is an optimal subspace for the decision variable $g$. The KFFL viewpoint also suggests that the innovation may play a useful role in exciting the system, which could potentially relax classical persistence-of-excitation requirements.








\bibliography{ifacconf}             

@article{Willems2005,
  title = {A note on persistency of excitation},
  volume = {54},
  ISSN = {0167-6911},
  number = {4},
  journal = {Systems and Control Letters},
  publisher = {Elsevier BV},
  author = {Willems,  Jan C. and Rapisarda,  Paolo and Markovsky,  Ivan and De Moor,  Bart L.M.},
  year = {2005},
  month = apr,
  pages = {325–329}
}

@inproceedings{Coulson2019,
  title = {Data-Enabled Predictive Control: In the Shallows of the DeePC},
  booktitle = {2019 18th European Control Conference (ECC)},
  publisher = {IEEE},
  author = {Coulson,  Jeremy and Lygeros,  John and Dörfler,  Florian},
  year = {2019},
  month = jun,
  pages = {307–312}
}

@article{innoop,
  title = {Data-driven output prediction and control of stochastic systems: An innovation-based approach},
  volume = {171},
  ISSN = {0005-1098},
  journal = {Automatica},
  publisher = {Elsevier BV},
  author = {Wang,  Yibo and You,  Keyou and Huang,  Dexian and Shang,  Chao},
  year = {2025},
  month = jan,
  pages = {111897}
}

@article{bridging,
  title = {Bridging Direct and Indirect Data-Driven Control Formulations via Regularizations and Relaxations},
  volume = {68},
  ISSN = {2334-3303},
  number = {2},
  journal = {IEEE Transactions on Automatic Control},
  publisher = {Institute of Electrical and Electronics Engineers (IEEE)},
  author = {D\"orfler,  Florian and Coulson,  Jeremy and Markovsky,  Ivan},
  year = {2023},
  month = feb,
  pages = {883–897}
}

@article{Mattsson2023,
  title = {On the regularization in {D}ee{PC}},
  volume = {56},
  ISSN = {2405-8963},
  number = {2},
  journal = {IFAC-PapersOnLine},
  publisher = {Elsevier BV},
  author = {Mattsson,  Per and Sch\"{o}n,  Thomas B.},
  year = {2023},
  pages = {625–631}
}

@article{Mattsson2024,
  title = {On the Equivalence of Direct and Indirect Data-Driven Predictive Control Approaches},
  volume = {8},
  ISSN = {2475-1456},
  journal = {IEEE Control Systems Letters},
  publisher = {Institute of Electrical and Electronics Engineers (IEEE)},
  author = {Mattsson,  Per and Bonassi,  Fabio and Breschi,  Valentina and Sch\"{o}n,  Thomas B.},
  year = {2024},
  pages = {796–801}
}

@article{SPC,
  title = {{SPC}: Subspace Predictive Control},
  volume = {32},
  ISSN = {1474-6670},
  number = {2},
  journal = {IFAC Proceedings Volumes},
  publisher = {Elsevier BV},
  author = {Favoreel,  Wouter and Moor,  Bart De and Gevers,  Michel},
  year = {1999},
  month = jul,
  pages = {4004–4009}
}

@inproceedings{Breschi2023cdc,
  title = {On the Impact of Regularization in Data-Driven Predictive Control},
  booktitle = {2023 62nd IEEE Conference on Decision and Control (CDC)},
  publisher = {IEEE},
  author = {Breschi,  Valentina and Chiuso,  Alessandro and Fabris,  Marco and Formentin,  Simone},
  year = {2023},
  month = dec,
  pages = {3061–3066}
}

@article{Breschi2023,
  title = {Data-driven predictive control in a stochastic setting: a unified framework},
  volume = {152},
  ISSN = {0005-1098},
  journal = {Automatica},
  publisher = {Elsevier BV},
  author = {Breschi,  Valentina and Chiuso,  Alessandro and Formentin,  Simone},
  year = {2023},
  month = jun,
  pages = {110961}
}

@inproceedings{IV1,
  title = {Data-enabled predictive control with instrumental variables: the direct equivalence with subspace predictive control},
  booktitle = {2022 IEEE 61st Conference on Decision and Control (CDC)},
  publisher = {IEEE},
  author = {van Wingerden,  Jan-Willem and Mulders,  Sebastiaan P. and Dinkla,  Rogier and Oomen,  Tom and Verhaegen,  Michel},
  year = {2022},
  month = dec,
  pages = {2111–2116}
}

@article{IV2,
  title = {Closed-loop Aspects of Data-Enabled Predictive Control},
  volume = {56},
  ISSN = {2405-8963},
  number = {2},
  journal = {IFAC-PapersOnLine},
  publisher = {Elsevier BV},
  author = {Dinkla,  Rogier and Mulders,  Sebastiaan P. and van Wingerden,  Jan-Willem and Oomen,  Tom},
  year = {2023},
  pages = {1388–1393}
}

@article{IV3Dinkla2026,
  title = {Closed-loop data-enabled predictive control and its equivalence with closed-loop subspace predictive control},
  volume = {183},
  ISSN = {0005-1098},
  journal = {Automatica},
  publisher = {Elsevier BV},
  author = {Dinkla,  Rogier and Oomen,  Tom and Mulders,  Sebastiaan Paul and van Wingerden,  Jan-Willem},
  year = {2026},
  month = jan,
  pages = {112556}
}

@article{IV_Wang2023,
  title = {Data-Driven Predictive Control Using Closed-Loop Data: An Instrumental Variable Approach},
  volume = {7},
  ISSN = {2475-1456},
  journal = {IEEE Control Systems Letters},
  publisher = {Institute of Electrical and Electronics Engineers (IEEE)},
  author = {Wang,  Yibo and Qiu,  Yiwen and Sader,  Malika and Huang,  Dexian and Shang,  Chao},
  year = {2023},
  pages = {3639–3644}
}

\appendix
\section{$\gamma$DDPC \& innovation nullspace} \label{sec:appA}   

We recall that, in the LQ factorization, $L$ is lower block triangular with invertible diagonal blocks $L_{ii}$, and $Q$ is orthogonal with block rows satisfying $Q_iQ_i^\top=I$, $Q_iQ_j^\top=0$ for $i\neq j$. We define the decision variables $\gamma = Qg$.

Using the first block row, $\gamma_1$ is fixed by $\text{col}(\mathbf{u}_{\text{ini}}, \mathbf{y}_{\text{ini}}) = L_{11}\gamma_1$. The second and third block rows give
\begin{equation} 
\begin{bmatrix}
            \mathbf{u} \\\hat{\mathbf{y}}
        \end{bmatrix}=\begin{bmatrix}
            L_{21} & L_{22} & 0 \\
            L_{31} & L_{32} & L_{33}
        \end{bmatrix} \begin{bmatrix}
            L_{11}^{-1} \begin{bmatrix}
                \mathbf{u}_{\text{ini}} \\\mathbf{y}_{\text{ini}}
            \end{bmatrix} \\ 
            \gamma_2 \\ \gamma_3
        \end{bmatrix}
\end{equation}
From this, the predictor decomposes as $\hat{\mathbf{y}} =\hat{\mathbf{y}}_{\text{LS}}+L_{33} \gamma_3$ in (\ref{eq:gamma_de}). Thus $\gamma_3$ parametrizes the deviation from SPC/OLS predictor.

Now we define $\Phi = \text{col}(U_p,Y_p,U_f)=M_1P_1$, with $M_1$ invertible and $P_1 = \text{col}(Q_1,Q_2)$. Then the orthogonal projector onto $\text{range}(\Phi)$ is $\Pi = P_1^\top P_1$, so $I-\Pi = Q_3^\top Q_3$. Using the third block row of the LQ factorization, we obtain $\hat{E}_{f,\text{LS}}  = Y_f(I-\Pi) = L_{33}Q_3$, and hence $L_{33}L_{33}^\top = Y_f(I-\Pi)Y_f^\top = \hat{E}_{f,\text{LS}}\hat{E}_{f,\text{LS}}^\top$.

Geometrically, $\gamma_3$ lives in the residual (innovation-like) subspace spanned by $\hat{E}_{f,\text{LS}}$. Penalizing $\Vert \gamma_3 \Vert^2$ in the $\gamma$-DDPC cost therefore softly suppresses these innovation directions, and restricts $g$ towards the approximate innovation null space $\text{null}(\hat{E}_{f,\text{LS}})$.

\end{document}